\documentclass{amsart}
\usepackage{amsmath}
\usepackage{amssymb}
\usepackage{amsopn}
\usepackage{epsfig}
\usepackage{amsfonts}
\usepackage{latexsym}
\usepackage{graphicx}
\newcommand{\prob}[1]{\ensuremath{{\rm P}_{\!\In}\left( #1 \right)}}

\newcommand{\expec}[1]{\ensuremath{{\rm E}_{\In}\!\left[#1\right]}}

\newcommand{\varfr}[1]{\ensuremath{{\rm Var}_{_{\scriptstyle\f}}\!\left( #1 \right)}}
\newcommand{\varab}[1]{\ensuremath{{\rm Var}_{_{\scriptstyle \ab}}\!\left( #1 \right)}}
\newcommand{\var}[1]{\ensuremath{{\rm Var}_{\In}\!\left( #1 \right)}}

\newcommand{\In}{\nu}       
\newcommand{\st}{\pi}    
\newcommand{\ei}{\gamma}   
\newcommand{\ab}{{\textsc{A}}}          
\newcommand{\f}{{\textsc{F}}}          
\newcommand{\ea}{\varepsilon_{\ab}}          
\newcommand{\ef}{\varepsilon_{\f}}          
\newcommand{\free}{free\;}          
\newcommand{\sorbed}{adsorbed\;}    
\newcommand{\thab}{\theta_\ab}
\newcommand{\thf}{\theta_\f}
\newcommand{\sab}{S^{\ab}}
\renewcommand{\sf}{S^{\f}}
\newcommand{\chab}{\varphi^{\ab}}
\newcommand{\chf}{\varphi^{\f}}
\newcommand{\jchab}{\hat{\varphi}^{\ab}}
\newcommand{\jchf}{\hat{\varphi}^{\f}}

\newcommand{\gf}{G^{\f}}

\newcommand{\fab}{f^{\ab}}
\newcommand{\ff}{f^{\f}}
\newcommand{\jab}{\hat{f}^{\ab}}
\newcommand{\jf}{\hat{f}^{\f}}

\newcommand{\MBDP}[1]{Proposition\,#1  in \cite{MBD}}
\newcommand{\MBDE}[1]{Equation\,(#1) in \cite{MBD}}

\newcommand{\A}{\mathcal{A}}

\newcommand{\La}{\Lambda}
\newcommand{\me}{\mathrm{e}}
\newcommand{\md}{\mathrm{d}}

\theoremstyle{plain}
\newtheorem{theorem}{Theorem}[section]

\newtheorem{lemma}{Lemma}[section]
\newtheorem{proposition}{Proposition}[section]

\title[]{A simple stochastic reactive transport model}
\author[]{Michel Dekking}
\author[]{Derong Kong}

\address{3TU Applied Mathematics Institute and Delft
University of Technology, Faculty EWI, P.O.~Box 5031, 2600 GA
Delft, The Netherlands.}\email{F.M.Dekking@tudelft.nl, \quad D.Kong@tudelft.nl.}

\thanks{The first author would like to thank J.Bruining, G.Uffink and C.Kraaikamp for numerous enlightening conversations. The second author is partially supported  by the National
Natural Science Foundation of China 10971069 and Shanghai Education Committee Project 11ZZ41.
}

\date{\today}

\begin{document}
\maketitle

\begin{abstract}
We introduce a discrete time microscopic single particle model for kinetic transport. The kinetics is modeled by a two-state Markov chain, the transport by deterministic advection plus a random space step. The position of the particle after $n$ time steps is  given by a random sum of space steps, where the size of the sum is given by a Markov binomial distribution (MBD). We prove that by letting the length of the time steps and the intensity of the switching between states tend to zero linearly, we obtain a random variable $S(t)$, which is closely connected to a well known (deterministic)
 PDE reactive transport model from
 the civil engineering literature. Our model explains (via bimodality of the MBD) the  double peaking behavior
  of the concentration of the free part of solutes in the PDE model. Moreover, we show for instantaneous injection of the solute that the partial densities of the  free and adsorbed part of the solute at time $t$ do exist, and  satisfy the partial differential equations.

\smallskip


\medskip


 \end{abstract}

\section{Introduction}

We consider a mathematical model for the displacement of a solute through a medium which apart from a constant flow (advection) and a dispersion (diffusion) interacts with the medium by intermittent adsorption (the kinetics). Our goal is to connect  a stochastic single particle model to the well known deterministic model which describes this process by a pair of partial differential equations.

In Section~\ref{sec: PDE} we give an introduction to the  deterministic reactive
transport model (as e.g. in \cite{Mich}) characterized by a pair of partial
differential equations.

 In Section~\ref{sec: SRTM} we give our simple discrete time microscopic single particle  stochastic reactive transport
 model.
In Section~\ref{sec: prop Kn} we calculate the probability
generating functions of the Markov binomial distribution (MBD) which is described in
Section~\ref{sec: SRTM}. These are helpful
 to consider the convergence of our simple discrete time stochastic model
 by letting the time step go to zero. This will be discussed in Section~\ref{sec: towards condtious time}.
In Section~\ref{sec:continuous time model} we compare our discrete time model with the obvious continuous time model.

In Section~\ref{sec: densPDE} we show for instantaneous injection of the solute that the partial probability densities of the free and adsorbed parts of the solute do
satisfy the PDE's defined in Section~\ref{sec: PDE}.
In Section~\ref{sec: mom} we compute the means and variances of our stochastic reactive transport model. Actually our formula fills a gap in
\cite{Mich}: since the authors erroneously state that the variances
are linear in the initial distribution, they only give the result
for two initial distributions (this might be connected to their
formula (22), which is incorrect).

In  Section~\ref{sec:Double Peaking} we study the probability density function of
our stochastic reactive transport model. This gives us a new and more precise  point of view
at the double peaking behavior in the concentration of the free part of the solute
discussed by  Michalak and Kitanidis in \cite{Mich}.

\section{The PDE reactive transport model}\label{sec: PDE}
\begin{figure}[h]
  \centering {\includegraphics[width=5cm]{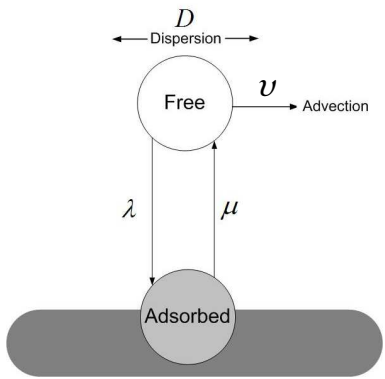}\\
  \caption{The schematic description of the kinetic transport model.\label{figure:01}}}
\end{figure}

We describe shortly the model used by Michalak and Kitanidis  in \cite{Mich} (see \cite{Lindstrom} for a more extensive treatment).
Given is a solute that has a sorbed part that does not move, and a free part that moves
in the $x$-direction by advection and  dispersion.
Let $C_\f(t,x)$ and $C_\ab(t,x)$ denote the concentration functions
 of the \free and the \sorbed part of the solute at time $t$ at position $x$.
By applying mass conservation and Fick's law one can set up the following
pair of differential equations:
 \begin{equation}\label{Dequation}
 \begin{split}
 \frac{{\partial C_\f(t,x)
}}{{\partial t}}+\frac{{\partial C_\ab(t,x) }}{{\partial t}}  &= D\frac{{\partial^2 C_\f(t,x)}} {{\partial x^2 }}-v\frac{{\partial C_\f(t,x)}} {{\partial x }}
 ,\\
\frac{{\partial C_\ab(t,x) }}{{\partial t}} &= - \mu C_\ab(t,x) +
\lambda C_\f(t,x).
\end{split}
\end{equation}
Here $D$ is called the dispersion coefficient and $v$ the advection velocity.
The parameters
$\lambda$ and $\mu$ denote the rates of changes as described in Figure \ref{figure:01}, with $\lambda$ for the change from
free to adsorbed and $\mu$ for the change from adsorbed to free.

The initial and boundary conditions are given by
\begin{equation*}
\begin{split}
C_\tau(0,x)=\In_\tau\delta(x),\quad
\lim_{x\rightarrow\infty}C_\tau(t,x)=\lim_{x\rightarrow\infty}\frac{\partial C_\tau(t,x)}{\partial x}=0\quad\textrm{for}~t\ge 0,~\tau\in\{\f,\ab\}
\end{split}
\end{equation*}
where $(\In_\f,\In_\ab)$ is a probability vector and $\delta$  the Dirac delta function.

 Michalak and Kitanidis have a slightly different set up, where the basic quantities are the aqueous concentration $C$ and the contaminant mass sorbed per mass of aquifer solids $S$. The connection is given by
$$C_\f=\eta C,\quad C_\ab=\rho S,$$
where $\eta$ is the porosity and $\rho$ mass of aquifer solids per total volume.

Also, Michalak and Kitanidis do not directly use $\lambda$ and $\mu$, but rather consider a distribution coefficient $K_d$ and a mass transfer coefficient $k$, which are given by
$$
\lambda=\frac{\rho K_d}{\eta}k,\quad \mu=k.
$$

The  main goal of the authors of \cite{Mich} is to obtain closed form expressions  for
the $m^{th}$ normalized moments for the \free and \sorbed phase, defined by
$$
 {M}_{\tau}^{(m)}(t)=\frac{1}{M_{\tau}^{(
0)}(t)}\int\limits_{-\infty
}^{+\infty}{x^{m}C_\tau(t,x)\,\mathrm{d}x},\quad\tau\in\{\f,\ab\},
$$
where the normalizing constants are given by
$
{M}_{\tau}^{(0)}(t)={\displaystyle\int\limits}C_\tau(t,x)\,\md
x.
$

These moments (for $m=1$ and $m=2$) are obtained in \cite{Mich} by taking Fourier transforms in the partial differential equations (\ref{Dequation}), and differentiating.
We copy here the formula\footnote{the $\oplus$ is a $+$ sign in (\cite{Mich}), but should be a $-$sign} from (\cite{Mich}, page 2136)
 for the normalized second central moment $\mu_2^*(t)$ where the solute is in the \free phase both at time 0 \emph{and} at time $t$:
\begin{equation}\label{eq:Michalak}
\begin{split}
\mu_2^*(t)&={\frac {{t}^{2}\A{v}^{2}\beta\, \left( \beta-1 \right) ^{2}}{ \left(
\beta+1 \right) ^{2} \left( 1+\beta\,\A \right) ^{2}}}+t \left( {
\frac {2{\it D}}{\beta+1}}+{\frac {2{v}^{2}\beta}{k \left( \beta+1
 \right) ^{3}}} \right)\\
 &\hspace{1cm}+ t \A \left( {\frac {4{v}^{2}\beta\, \left( -{
\beta}^{2}\A-{\beta}^{2}-\beta+1 \right) }{k \left( 1+\beta\,\A \right)
^{2} \left( \beta+1 \right) ^{3}}}+{\frac {2{\it D}\,\beta\,
 \left( \beta-1 \right) }{ \left( \beta+1 \right)  \left( 1+\beta\,\A
 \right) }} \right) \\
 &\hspace{1.5cm}+ {\frac {2{v}^{2}\beta\, \left( 1-\A \right)
 \left( 3\,{\beta}^{2}\A-3-\beta\, ( \A\oplus 1 )  \right) }{{k}^{2
} \left( 1+\beta\,\A \right) ^{2} \left( \beta+1 \right) ^{4}}}+{
\frac {4{\it D}\,\beta\, \left( 1-\A \right) }{k \left( 1+\beta\,\A
 \right)  \left( \beta+1 \right) ^{2}}}.
 \end{split}
\end{equation}
Here Michalak and Kitanidis have made the following abbreviations:
$$\beta=\frac{\rho K_d}{\eta}=\frac{\lambda}{\mu},\quad \A=\A(t)=\exp(-(\beta+1)kt)=\exp(-(\lambda+\mu)t).$$

\section{A simple stochastic reactive transport model}\label{sec: SRTM}

We describe the behavior of a single particle in the solute. Time $t$ is
 discretized  by choosing some $n$, and dividing $[0,t]$ into
$n$ intervals of the same length $$\Delta t=t/n.$$
  We suppose in such an interval of length $\Delta t$ that the
particle can only be in one of the two states: `\free' or `\sorbed',
which we code by the letters $\f$ and $\ab$. The particle can only
move when it is `\free', and in this case its displacement has two
components: dispersion and advection.
Let $X_k, k\ge 1$ be the displacement of the particle due to the dispersion
the $k${th} time that it is `free'. We model the $X_k$ as independent
identically distributed random variables satisfying
\begin{equation}\label{eq:displace}
\expec{X_k}=0, \quad\expec{X_k^2}=2 D\Delta t\quad {\rm and } \quad   \expec{X_k^3}=o(\Delta t)\;\mathrm{as}\; \Delta t \downarrow 0,
\end{equation}
where $D>0$, and $\In=(\In_\f,\In_\ab)$ is the initial distribution describing the state of the particle at time 0.
 When the
particle is free during the interval $[(k-1)\Delta t, k\Delta t]$
for some $k$, the displacement due to advection is given by $v\Delta
t$ with $v$ the (deterministic) advection velocity.

In order to model the kinetics,
let $\{Y_k, k\ge 1\}$ be a process taking values in $\{\f,\ab\}$ (we will make a choice for $\{Y_k\}$ below), and let
  $$K_n=\sum_{k=1}^n {\bf 1}_{\{Y_k=\f\}}$$
   be the occupation time of the process $\{Y_k\}$  in state $F$ up to time $n$.
\begin{figure}[h]
  \centering {\includegraphics[width=6.5cm]{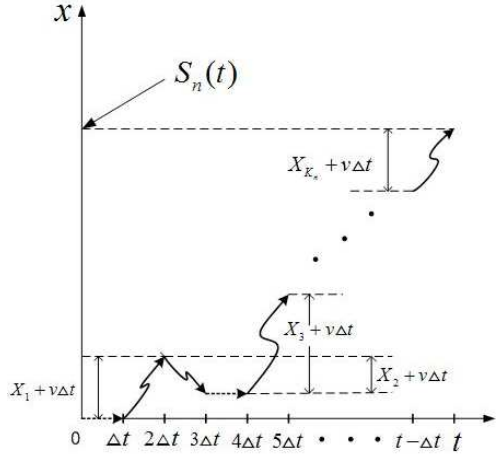}\\
  \caption{The position $S_n(t)$ of the particle at time $t=n\Delta t$ with $Y_1=\ab, Y_2=\f, Y_3=\f, Y_4=\ab, Y_5=\f,\cdots, Y_n=\f.$\label{figure:02}}}
\end{figure}

   Now let $S_n(t)$ be the position of the particle at time $t=n \Delta t$. Then by the above (see also Figure \ref{figure:02})
we can write $S_n(t)$ as
\begin{equation*}\label{S_n}
S_n(t)=\sum_{k=1}^{K_n} (X_k +v\Delta t).
\end{equation*}
Here we assume that $K_n$ is independent of the dispersion $X_k,\, k=1,\dots,K_n$.

We want to compare our stochastic model with the PDE-model of Michalak and Kitanidis from Section~\ref{sec: PDE}. Since these authors consider the solute with given states
 (`free' or `adsorbed')  at time $t$, we need to consider the conditional random variables $S^\f_n(t)$ and $S_n^\ab(t)$, i.e., the position of the particle at time $t$ given that it is `free' and `adsorbed' respectively at time $t=n\Delta t$. Let $K_n^\tau$ be the random variable $K_n$ conditioned on $Y_n=\tau$ with $\tau\in\{\f,\ab\}$, i.e., $K_n^\tau$ counts the number of intervals $[(k-1)\Delta t, k\Delta t], 1\le k\le n$ where the particle is free, conditioned on the particle being in state $\tau$ in $[t-\Delta t, t]$.
 Then $S_n^\tau(t)$ can be written as
 $$
S_n^\tau(t)=\sum_{k=1}^{K_n^\tau}(X_k+v\Delta t).
 $$

 The distributions of $K_n$ and $K_n^\tau$ are determined by the process $\{Y_k\}$.
We  take for  $\{Y_k, k\ge 1\}$ a Markov chain  on the two states $\{\f,\ab\}$ with initial distribution $\In=(\In_\f, \In_\ab)$ and transition matrix
 \begin{equation}\label{transition matrix}
P=\left[ {\begin{array}{*{20}c}
   {P(\f,\f) } & {P(\f,\ab) }  \\
   {P(\ab,\f) } & {P(\ab,\ab) }  \\
\end{array}}  \right] =\left[ {\begin{array}{*{20}c}
   {1 - a} & {a}  \\
   {b} & {1-b}  \\
\end{array}} \right],
\end{equation}
where we assume $0<a, b<1$. The distribution of $K_n$ is then well
known, and is called a Markov binomial distribution (MBD)
 (see, e.g., \cite{MBD,Omey}).

Clearly the stationary distribution   $(\st_\f,\st_\ab)$  of the
Markov chain $\{Y_k,k\ge 1\}$ is given by
$$\st_\f=\frac{b}{a+b}, \quad \st_\ab=\frac{a}{a+b}.$$
It is useful to consider the  \emph{excentricities} $\ef$ and
$\ea$ of an initial distribution $\In$ given by
$$ \varepsilon_\tau:=\varepsilon_\tau(\In)=1-\frac{\In_\tau}{\st_\tau},\quad\mbox{for}\quad \tau\in\{\f, \ab\}.$$
We can then write for $k\ge 1$
\begin{equation*}\label{P(Y_k)}
\prob{Y_k=\tau}=\st_\tau(1-\varepsilon_\tau\,\ei^{k-1}),
\end{equation*}
where $\gamma=1-a-b$ is the smallest eigenvalue of $P$ (see also \cite{MBD} for the computations).

\section{Probability generating functions of $K_n$ and $K_n^\tau$}\label{sec: prop Kn}
We compute in this section the probability generating functions of
$K_n$ and $K_n^\tau$. These  are useful when we consider the convergence
of the random variables $S_n(t)$ and $S_n^\tau(t)$ as $n$ goes to infinity.

 Given $n\ge 1$, let $f_{n}$ be the probability mass function of $K_n$, i.e.,
 $$f_{n}(j)=\prob{K_n=j}. $$
In particular $f_n(j)=0$ if $j<0$ or $j>n$.
Straightforward computations as in \cite{Viveros} or \cite{MBD} yield that
\begin{equation*}\label{eq:recursion_prob_mass_func}
f_{n+2}(j+1)=(1-b)f_{{n+1}}(j+1)+(1-a)f_{{n+1}}(j)-(1-a-b)f_{{n}}(j)
\end{equation*}
with initial conditions
\begin{equation*}\label{eq:initial conditions}
\begin{split}
&f_1(0)=\In_\ab,\quad f_{1}(1)=\In_\f;\\
&f_2(0)=\In_\ab (1-b),\quad f_2(1)=\In_\ab b+\In_\f a,\quad f_2(2)=\In_\f(1-a).
\end{split}
\end{equation*}
Let $G_n$ be the probability generating function of $K_n$, i.e.,
$$
G_n(s)=\expec{s^{K_n}}=\sum_{j=0}^n f_n(j)s^j.
$$
It follows from the above recursion equation for $f_n$  that
$$
G_{n+2}(s)=\big((1-a)s+(1-b)\big)G_{n+1}(s)-(1-a-b)s G_n(s)
$$
with initial conditions
$$
G_1(s)=\In_\ab+\In_\f s,\quad G_2(s)=\In_\ab(1-b)+(\In_\ab b+\In_\f a)s+\In_\f(1-a)s^2.
$$
By solving the difference equation of $G_n$ with the initial conditions we obtain the probability generating function of $K_n$
(see also \cite{ Viveros}).
\begin{equation}\label{eq:generating func}
  \begin{split}
   G_n(s)=&~\frac{\In_\ab\big(1-\beta(s)+b(s-1)\big)+\In_\f
   s\big(a-\beta(s)+s(1-a)\big)}{\alpha(s)-\beta(s)}\,\alpha(s)^{n-1}\\
   &+\frac{\In_\ab\big(1-\alpha(s)+b(s-1)\big)+\In_\f
   s\big(a-\alpha(s)+s(1-a)\big)}{\beta(s)-\alpha(s)}\,\beta(s)^{n-1},
  \end{split}
\end{equation}
where
\begin{equation}\label{eq:alpha_beta}
\begin{split}
\alpha(s)&=\frac{1}{2}\Big((1-a)s+(1-b)+\sqrt{\big((1-a)s-(1-b)\big)^2+4
a b
s}\Big),\\
\beta(s)&=\frac{1}{2}\Big((1-a)s+(1-b)-\sqrt{\big((1-a)s-(1-b)\big)^2+4
a b s}\Big).
\end{split}
\end{equation}

Next we are going to consider the probability generating function of $K_n^\tau$ for $\tau\in\{\f,\ab\}$.
Given $n\ge 1$, let $f^{\tau}_n$ be the  probability mass function of $K_n^\tau$, i.e.,
\begin{equation}\label{eq:fntau}
f_n^{\tau}(j)=\prob{K_n^\tau=j}=\prob{K_n=j~|~ Y_n=\tau}.
\end{equation}
In order to deal with $f_n^\tau$ it is simpler to deal with the \emph{partial probability mass functions}
$$
\hat{f}_n^\tau(j)=\prob{K_n=j, Y_n=\tau}=f_n^\tau(j) \prob{Y_n=\tau},
$$
since these satisfy the \emph{same} recursion equation  as $f_n$:
\begin{equation*}\label{eq:conditional recursion equation}
\hat{f}_{n+2}^{\tau}(j+1)=(1-b)\hat{f}_{n+1}^{\tau}(j+1)+(1-a)\hat{f}_{n+1}^{\tau}(j)-(1-a-b)
\hat{f}_{n}^{\tau}(j).
\end{equation*}
Only the initial conditions are different:
\begin{equation*}\label{eq:inital_F}
  \begin{split}
    \jf_1(0)&=0,\quad \jf_1(1)=\In_\f;\\
    \jf_2(0)&=0, \quad\jf_2(1)=\In_\ab b, \quad\jf_2(2)=\In_\f(1-a);
  \end{split}
\end{equation*}
and
\begin{equation*}\label{eq:inital_A}
  \begin{split}
    \jab_1(0)&=\In_\ab, \quad\jab_1(1)=0;\\
    \jab_2(0)&=\In_\ab(1-b), \quad\jab_2(1)=\In_\f a,\quad \jab_2(2)=0.
  \end{split}
\end{equation*}
Then using the above recursion equation of $\hat{f}_n^\tau$ with these initial conditions,  the probability generating function $G_n^\tau$ of $K_n^\tau$ can be obtained in a similar way as for $G_n$ (see also \cite{Viveros}).
\begin{equation}\label{eq:generating func_free}
\begin{split}
   &G_n^\f(s)=\sum_{j=0}^n \ff_n(j)s^j=\sum_{j=0}^n \frac{\jf_n(j)}{\prob{Y_n=\f}}s^j=\frac{1}{\st_\f(1-\ef\ei^{n-1})}\sum_{j=0}^n \jf_n(j)s^j\\
   &=\frac{\In_\ab b-\In_\f\beta(s)+\In_\f(1-a)s}{(\alpha(s)-\beta(s))\st_\f(1-\ef\ei^{n-1})}\,s\alpha(s)^{n-1}+\frac{\In_\ab b-\In_\f \alpha(s)+\In_\f(1-a)s}{(\beta(s)-\alpha(s))\st_\f(1-\ef\ei^{n-1})}\,s\beta(s)^{n-1},
   \end{split}
\end{equation}
and
\begin{equation}\label{eq:generating func_absorbed}
  \begin{split}
   &G_n^\ab(s)=\sum_{j=0}^n \fab_n(j)s^j=\sum_{j=0}^n \frac{\jab_n(j)}{\prob{Y_n=\ab}}s^j=\frac{1}{\st_\ab(1-\ea\ei^{n-1})}\sum_{j=0}^n \jab_n(j)s^j\\
   &=\frac{\In_\ab(1- b)+\In_\f a s-\In_\ab\beta(s)}{(\alpha(s)-\beta(s))\st_\ab(1-\ea\ei^{n-1})}\,\alpha(s)^{n-1}+\frac{\In_\ab (1-b)+\In_\f a s -\In_\ab\alpha(s)}{(\beta(s)-\alpha(s))\st_\ab(1-\ea\ei^{n-1})}\,\beta(s)^{n-1}.
   \end{split}
\end{equation}

\section{Towards continuous time}\label{sec: towards condtious time}

 To get closer to the PDE model in  Section~\ref{sec: PDE}, we have to fix $t=n\Delta t>0$ and then let the time step $\Delta t$ tend to
 $0$. Hence $n$ goes to infinity. We consider the rates of changes  $\lambda$ and $\mu$ from Section \ref{sec: PDE}.
 Since the probability that a particle changes its state is
proportional to the length of the time step $\Delta t$ (if $\Delta
t$ is small), we should put
\begin{equation}\label{Delta t}
a=\lambda \Delta t=\frac{\lambda t}{n},\qquad b=\mu \Delta
t=\frac{\mu t}{n}
\end{equation}
in the transition matrix $P$ in (\ref{transition matrix}).
Under this assumption, we will show in this
 section that the random variables $S_n(t)$ and $S_n^\tau(t)$ defined in Section \ref{sec:
 SRTM} converge in distribution to some random variables $S(t)$ and
 $S^\tau(t)$ respectively. To achieve this,
we first consider the characteristic function $\varphi_{t,n}$ of $S_n(t)$,
i.e.,
$$
\varphi_{t,n}(u)=\expec{\me^{i u S_n(t)}}=\expec{\me^{i u \sum_{k=1}^{K_n}
(X_k+v\Delta t)}}=G_n\Big(\expec{\me^{iu(X_1+v\Delta t)}}\Big),
$$
where $G_n$ is the generating function of $K_n$ given in
(\ref{eq:generating func}). It is  well known that
(cf.~\cite{Billing})
\begin{equation*}
  \begin{split}
   &\Big| \expec{\me^{i u(X_1+v \Delta t)}}-\sum_{k=0}^2\frac{(i u)^k}{k!}\expec{(X_1+v\Delta t)^k}\Big|\\
   \le &\; \frac{|u|^3}{3!}\expec{\big|X_1+v\Delta t\big|^3}   =  o(\Delta t)=o\big(\frac{1}{n}\big),
  \end{split}
\end{equation*}
where the last equalities hold since $\expec{X_1^3}=o(\Delta t)$ in Equation (\ref{eq:displace}) and $t$ is always assumed to be fixed.  We then obtain by (\ref{eq:displace}) that
\begin{equation}\label{eq:approx X1+vt/n}
\begin{split}
 \expec{\me^{iu(X_1+v\Delta t)}}&=1+i u\expec{X_1+v\Delta
t}-\frac{u^2}{2}\expec{(X_1+v\Delta
t)^2}+o\big(\frac{1}{n}\big)\\
&=1+\frac{t u(i v -D  u)}{n}+o\big(\frac{1}{n}\big)=1+\frac{z}{n}+o\big(\frac{1}{n}\big),
\end{split}
\end{equation}
defining
$$ z:=t u(i v -D  u).$$

 Substituting (\ref{Delta t}) and (\ref{eq:approx X1+vt/n}) into Equation (\ref{eq:alpha_beta}) and letting $n$ go to infinity, we obtain
{\small\begin{equation*}
\begin{split}
  &\quad\alpha\Big(\expec{\me^{i u(X_1+v\Delta t)}}\Big)^n=\alpha\Big(1+\frac{z}{n}+o\big(\frac{1}{n}\big)\Big)^{\!n}\\
  &=\left(\frac{(1-\frac{\lambda t}{n})(1+\frac{z}{n})+(1-\frac{\mu t}{n})+\sqrt{\big((1-\frac{\lambda t}{n})(1+\frac{z}{n})-(1-\frac{\mu t}{n})\big)^2+\frac{4\lambda\mu t^2}{n^2}(1+\frac{z}{n})}}{2}+o\big(\frac{1}{n}\big)\right)^n\\
  &\longrightarrow\; \exp{\Big(\frac12\Big[ z-(\lambda+\mu)t+\sqrt{\big(z+t(\mu-\lambda)\big)^2+4\lambda\mu t^2}\Big]\Big)} \\
  &\quad =\exp{\Big(-\frac{t}{2}\Big[D u^2-i v u+\lambda+\mu-\sqrt{(Du^2-ivu+\lambda-\mu)^2+4\lambda\mu}\,\Big]\Big)}.
\end{split}
\end{equation*}}
Here we chose the complex square root of $(Du^2-ivu+\lambda-\mu)^2+4\lambda\mu$ with positive real part.

Similarly, the corresponding limit for $\beta\big(\expec{\me^{i u(X_1+v\Delta t)}}\big)^n$
  is obtained by replacing the minus in front of the square root of the last equality by a plus. It seems convenient to introduce the following two notations:
 \begin{equation}\label{thab thf}
\begin{split}
\thab=\thab(u):=&\sqrt{(Du^2-ivu+\lambda-\mu)^2+4\lambda\mu}-(D u^2-i v u),\\
\thf=\thf(u):=&\sqrt{(Du^2-ivu+\lambda-\mu)^2+4\lambda\mu}+(D u^2-i v u).
\end{split}
\end{equation}
Then the limits for $\alpha\big(\expec{\me^{i u(X_1+v\Delta t)}}\big)^n$ and $\beta\big(\expec{\me^{i u(X_1+v\Delta t)}}\big)^n$ can be rewritten as:
\begin{equation*}
\begin{split}
\lim\limits_{n\rightarrow\infty}\alpha\big(\expec{\me^{i u(X_1+v\Delta t)}}\big)^n&=\exp{\big((\thab-\lambda-\mu)t/2\big)},\\
\lim\limits_{n\rightarrow\infty}\beta\big(\expec{\me^{i u(X_1+v\Delta t)}}\big)^n&=\exp{\big(-(\thf+\lambda+\mu)t/2\big)}.
\end{split}
\end{equation*}
So we obtain by substituting (\ref{Delta t}) and (\ref{eq:approx X1+vt/n}) into Equation (\ref{eq:generating func}) that the limit of the characteristic functions $\varphi_{t,n}$
of $S_n(t)$ is  a function $\varphi_t$ given by

\begin{equation}\label{eq:charac func S(t)}
\varphi_t(u)=\me^{(\thab-\lambda-\mu)t/2}\frac{\In_\ab \thf+\In_\f\,\thab+\lambda+\mu}{\thab+\thf}+\me^{-(\thf+\lambda+\mu)t/2}\frac{\In_\ab\thab+\In_\f\,\thf-\lambda-\mu}{\thab+\thf}.
\end{equation}
 It is easy
to see that $\varphi_t$ is continuous at $u=0$. This implies that there exists a random variable, which we call $S(t)$, such that as $n\rightarrow\infty$
$$
S_n(t)\longrightarrow S(t)\quad\textrm{in~~distribution.}
$$

Next we are going to consider the convergence of the random variable
$\sf_n(t)$ as $n$ goes to infinity. In a similar way as for $S_n(t)$
we consider the characteristic function $\chf_{t,n}$ of $\sf_n(t)$, i.e.,
$$
\chf_{t,n}(u)=\expec{\me^{i u\sf_n(t)}}=\gf_n\Big(\expec{\me^{i u(X_1+v\Delta t)}}\Big),
$$
where $\gf_n$ is the  probability generating function of $K_n^\f$
given in Equation (\ref{eq:generating func_free}).
Substituting (\ref{Delta t})  and (\ref{eq:approx X1+vt/n}) into Equation (\ref{eq:generating func_free}) and letting $n$ go to infinity, we obtain that the limit of the characteristic function $\chf_{t,n}$ of $\sf_n(t)$ is a function $\chf_t$ given by
 \begin{equation}\label{eq:charac func S(t)_F}
   \chf_t(u)=\me^{(\thab-\lambda-\mu)t/2}\frac{\In_\f(\thab-\lambda-\mu)+2\mu}{(\thab+\thf)\st_\f(1-\ef\A)}
   +\me^{-(\thf+\lambda+\mu)t/2}\frac{\In_\f(\thf+\lambda+\mu)-2\mu}{(\thab+\thf)\st_\f(1-\ef\A)},
 \end{equation}
 where $\A=\A(t)=\exp{(-(\lambda+\mu)t)}$. Here we point out that the stationary distribution $(\st_\f,\st_\ab)$ and the excentricities $\ef,\ea$ do not depend on the time step $\Delta t$, since by (\ref{Delta t})
 $$
 \st_\f=\frac{b}{a+b}=\frac{\mu}{\lambda+\mu},~ \st_\ab=\frac{a}{a+b}=\frac{\lambda}{\lambda+\mu},~ \ef=1-\frac{\lambda+\mu}{\mu}\In_\f, ~ \ea=1-\frac{\lambda+\mu}{\lambda}\In_\ab.
 $$
  Again there exists a random variable, which we call $\sf(t)$, such that as $n\rightarrow\infty$
$$\sf_n(t)\longrightarrow \sf(t)\quad\textrm{in~~distribution.}$$

Similarly, substituting (\ref{Delta t})  and (\ref{eq:approx X1+vt/n}) into
Equation (\ref{eq:generating func_absorbed}) and letting $n$ go to
infinity, one can show  that there exists a random variable $\sab(t)$ such that $\sab_n(t)\longrightarrow \sab(t)$ in distribution as $n\rightarrow\infty$,
where $\sab(t)$ has
characteristic function:
 \begin{equation}\label{eq:charac func S(t)_A}
   \chab_t(u)=\me^{(\thab-\lambda-\mu)t/2}\frac{\In_\ab(\thf-\lambda-\mu)+2\lambda}{(\thab+\thf)\st_\ab(1-\ea\A)}
   +\me^{-(\thf+\lambda+\mu)t/2}\frac{\In_\ab(\thab+\lambda+\mu)-2\lambda}{(\thab+\thf)\st_\ab(1-\ea\A)}.
 \end{equation}

\section{Modeling the kinetics with a continuous time Markov chain}\label{sec:continuous time model}

In our model we used a simple discrete time set up.  This will be useful in Section \ref{sec:Double Peaking}, but it is worthwhile to compare our results with a model that involves a continuous time Markov chain. Let $Y(t), t\ge 0$ denote the state of the particle at time $t$.
Recall from Section \ref{sec: PDE}  that $\lambda$ and $\mu$ are the  rates of changes from `free' to `adsorbed' and `adsorbed' to `free' respectively. Hence it is natural to model the kinetics  by a two-state continuous time Markov chain $\{Y(t), t\ge 0\}$ with initial distribution $\prob{Y(0)=\tau}=\In_\tau, \tau\in\{\f,\ab\}$ and  generator matrix
$$
Q=\left(\begin{array}{cc}
  Q(\f,\f)&Q(\f,\ab)\\
  Q(\ab,\f)&Q(\ab,\ab)
\end{array}
\right)=\left(\begin{array}{cc}
  -\lambda&\lambda\\
  \mu&-\mu
\end{array}
\right).
$$
The solute can only move when it is free, and in this case we model the displacement due to dispersion and advection as a Brownian motion with drift $v$.

 A trick to deal with continuous time Markov chains  is  \emph{uniformization}. This idea gives us an alternative way to  model the $S(t)$ and $S^\tau(t), \tau\in\{\f,\ab\}$ obtained in Section \ref{sec: towards condtious time}.
Let $\La\ge\max{(\lambda,\mu)}$ be the rate of the uniformization. It follows that (see e.g.~\cite{Ross}, page 402) the continuous time Markov chain $\{Y(t), t\ge 0\}$ can be viewed as a discrete time Markov chain $\{Z_k, k\ge 0\}$ over the same state space $\{\f,\ab\}$ and the same initial distribution $\prob{Z_0=\tau}=\In_\tau, \tau\in\{\f,\ab\}$, but with the  transition matrix
$$
P_\La=\left(\begin{array}{cc}
  P(\f,\f)&P(\f,\ab)\\
  P(\ab,\f)&P(\ab,\ab)
\end{array}
\right)=\left(\begin{array}{cc}
  1-\lambda/\La&\lambda/\La\\
  \mu/\La&1-\mu/\La
\end{array}
\right).
$$
Let $N(t)$ be the number of the state transitions up to time $t$, which is a Poisson process with rate $\La$.  Let $K_{N(t)}$ be the occupation time of the chain $\{Z_k\}$ in state $\f$ up to time $N(t)$, which  is a Markov binomial distributed random variable, when conditioned on $N(t)$.
Since the solute can only move when it is free and the displacement in the free state is due to dispersion and advection, we model $\mathcal{X}_k$, the displacement during the $k${th} free interval, as a Brownian motion with drift $v$  stopped at time $T$ which is exponentially $\La$ distributed. So we put
 $$
\mathcal{ X}_k\stackrel{d}{=}\mathcal{N}(v T, {2 D}T)\quad\textrm{with}\quad T\stackrel{d}{=}Exp(\La).
 $$
Then we can write $\mathcal{H}_\La(t)$, the position of the particle at time $t$ with respect to the uniformization at rate $\La$,
as:
$$
\mathcal{H}_\La(t)=\sum_{k=1}^{K_{N(t)}}\mathcal{X}_k.
$$
Similarly, for $\tau\in\{\f,\ab\}$ we can define $\mathcal{H}_\La^\tau(t)$ merely by changing $K_{N(t)}$ to the conditional random variable $K_{N(t)}^\tau=K_{N(t)}~|~{\{Z_{N(t)}=\tau\}}$, i.e., $\mathcal{H}_\La^\tau(t)$ denotes the position of the particle at time $t$ conditioned on being in state $\tau$ at time $N(t)$.

Letting $\La$ go to infinity, one can show by using the characteristic functions of $\mathcal{H}_\La(t)$ and $ \mathcal{H}_\La^\tau(t)$ as in Section \ref{sec: towards condtious time} that
$$
\mathcal{H}_\La(t)\longrightarrow S(t),\quad \mathcal{H}_\La^\tau(t)\longrightarrow S^\tau(t)\quad\textrm{in~distribution},
$$
where $S(t), S^\tau(t)$ are the \emph{same} random variables as in Section \ref{sec: towards condtious time}.

 It is  even more natural  to look at the continuous time Markov chain $\{Y(t), t\ge 0\}$ directly.  Let $U(t)$ be the occupation time of the chain $\{Y(t)\}$ in state $\f$ up to time $t$, and let $f_{U(t)}$ be its probability density function. We model the displacement of the solute in the free phase as a Brownian motion with drift $v$. Then the position $H(t)$ of the particle at time $t$ can be written as a normal distribution with mean $v \, U(t)$ and variance $2D\, U(t)$.
Conditional on $U(t)$ it follows  from Equation (5) of \cite{Pedler} and Equation (\ref{eq:charac func S(t)}) that for $\phi>0$
\begin{equation*}
\begin{split}
    \int_0^\infty\expec{\me^{i u H(t)}}\me^{-\phi t}\md t
    =&~\int_0^\infty\int_0^\infty\me^{-(Du^2-i v u)x}\me^{-\phi t}f_{U(t)}(x)\,\md x\,\md t\\
    =&~(\In_\f,\In_\ab)\left(\begin{array}{cc}
      \lambda+\phi+Du^2-ivu&-\lambda\\
      -\mu&\mu+\phi
    \end{array}\right)^{-1}{1\choose 1}\\
    =&~\frac{\phi+\lambda+\mu+\In_\ab(Du^2-i v u)}{(\phi+\lambda+D u^2-i v u)(\phi+\mu)-\lambda\mu}=\int_0^\infty\expec{\me^{i u S(t)}} e^{-\phi t}\md t.
    \end{split}
\end{equation*}
Since  $\expec{\me^{i u S(t)}}$ is a continuous function of $t$ by Equation (\ref{eq:charac func S(t)}), it follows from Lerch's theorem (cf.~\cite{Schiff}, page 24) that $\expec{\me^{iu H(t)}}=\expec{\me^{iuS(t)}}$ for all $t\ge 0$. Hence $H(t)$ and $S(t)$ have the same distribution.

 Similarly, for $\tau\in\{\f,\ab\}$ let $H^\tau(t)=H(t)~|~{\{Y(t)=\tau\}}$ be the conditional random variable denoting  the position of the particle at time $t$ conditioned on being in  state  $\tau$ at time $t$.  From the proof of Theorem 1 in  \cite{Darroch}  one obtains that for $\phi>0$
\begin{equation*}
  \begin{split}
   \int_0^\infty\prob{Y(t)=\f}&\expec{\me^{iuH^\f(t)}}\me^{-\phi t}\md t   =\int_0^\infty\expec{\me^{iuH(t)}\textbf{1}_{\{Y(t)=\f\}}}\me^{-\phi t}\md t\\
   =&~(\In_\f,\In_\ab)\left(\begin{array}{cc}
      \lambda+\phi+Du^2-ivu&-\lambda\\
      -\mu&\mu+\phi
    \end{array}\right)^{-1}{1\choose 0}\\
   =&~\frac{\mu+\In_\f\phi}{(\phi+\lambda+D u^2-ivu)(\phi+\mu)-\lambda\mu}\\
   =&~\int_0^\infty\prob{Y(t)=\f}\expec{\me^{iu{S}^\f(t)}}\me^{-\phi t}\md t,
  \end{split}
\end{equation*}
where the last equality follows using Equation (\ref{eq:charac func S(t)_F}) and since  the limiting probability of a particle being in state $\tau$ at time $t$ is given by
\begin{equation}\label{eq:limit probability at time t}
\prob{Y(t)=\tau}=\lim_{n\rightarrow\infty}\prob{Y_n=\tau}=\lim_{n\rightarrow\infty}\st_\tau(1-\varepsilon_\tau\ei^{n-1})=\st_\tau(1-\varepsilon_\tau\A(t))
\end{equation}
with $\A(t)=\exp{(-(\lambda+\mu)t)}$.
Similarly one can also show that for $\phi>0$
\begin{equation*}
\begin{split}
\int_0^\infty\prob{Y(t)=\ab}\expec{\me^{iuH^\ab(t)}}\me^{-\phi t}\md t&=\frac{\lambda+\In_\ab(\phi+Du^2-ivu)}{(\phi+\lambda+D u^2-ivu)(\phi+\mu)-\lambda\mu}\\
&=\int_0^\infty\prob{Y(t)=\ab}\expec{\me^{iu{S}^\ab(t)}}\me^{-\phi t}\md t.
\end{split}
\end{equation*}
Again, using Lerch's theorem, it follows that $H^\tau(t)$ and ${S}^\tau(t)$ have the same distribution.
Therefore our discrete time model converges in distribution to the same random variables  as obtained by the natural continuous time Markov chain.

\section{Densities and partial differential equations}\label{sec: densPDE}
We will show in this section that for instantaneous injection of the solute, i.e., with
initial distribution $\In=(1,0)$, the partial probability density
functions $\jf_S(t,x)$ and $\jab_S(t,x)$ of $\sf(t)$ and $\sab(t)$
 do satisfy the
partial differential equations  in (\ref{Dequation}).

Let $f_S(t,x)$ and $f_S^\tau(t,x)$ denote respectively the probability density functions of $S(t)$ and $S^\tau(t)$ for $\tau\in\{\f,\ab\}$. Recall from (\ref{eq:limit probability at time t}) that the probability of a particle being in state $\tau$ at time $t$ is given by $\prob{Y(t)=\tau}=\st_\tau(1-\varepsilon_\tau\A(t))$. We define the partial probability density functions of $S^\tau(t)$ as
\begin{equation}\label{eq:partial densities}
  \hat{f}_S^\tau(t,x)=\prob{Y(t)=\tau}f^\tau_S(t,x)=\st_\tau(1-\varepsilon_\tau\A(t))f_S^\tau(t,x).
\end{equation}
Obviously $f_S(t,x)=\jf_S(t,x)+\jab_S(t,x)$.

\begin{lemma}\label{lemma:speed_thab}
Let $\thab=\thab(u),~ \thf=\thf(u)$ be defined as in {\rm(\ref{thab thf})}. Then
  $$\lim\limits_{u\rightarrow \infty}\big|u^2\big(\thab(u)-(\lambda-\mu)\big)\big|=2\lambda\mu/D,\quad\lim\limits_{u\rightarrow\infty}|\thf(u)/u^2|=2D.$$
\end{lemma}
\begin{proof}
 It is straightforward to check these formulas.
\end{proof}

\begin{lemma}\label{lemma:integrable_chf}
The probability  density function $\ff_S(t,\cdot)$ of $\sf(t)$ can be written as
 $$\ff_S(t,x)=\frac{1}{2\pi}\int \me^{-iux}\chf_t(u)\,\md u,$$
 where $\chf_t$ is the characteristic function of $\sf(t)$ given in {\rm(\ref{eq:charac func S(t)_F})}.
\end{lemma}
\begin{proof}
  We only need to show that $\chf_t$ is integrable. Obviously
  $\chf_t$ is a continuous function.
 So it suffices to show that
  $\int_{|u|>M}|\chf_t(u)|\,\md u<\infty$ for some $M>0$.
  From Lemma \ref{lemma:speed_thab} and (\ref{eq:charac func S(t)_F})  it follows that  for all $|u|$ large
  \begin{equation*}
  \begin{split}
      |\chf_t(u)|&\le\Big|\me^{(\thab-\lambda-\mu)t/2}\frac{\In_\f(\thab-\lambda-\mu)+2\mu}{(\thab+\thf)\st_\f(1-\ef\A)}\Big|
   +\Big|\me^{-(\thf+\lambda+\mu)t/2}\frac{\In_\f(\thf+\lambda+\mu)-2\mu}{(\thab+\thf)\st_\f(1-\ef\A)}\Big|\\
   &\le \frac{C_1}{u^2}+C_2 \me^{-D t  u^2/2},
   \end{split}
  \end{equation*}
 where $C_1, C_2$ are constants independent of $u$. This finishes the proof of the lemma.
\end{proof}
Surprisingly, Lemma \ref{lemma:integrable_chf} does not hold for $\sab(t)$,
but we still have the following.
\begin{lemma}\label{lemma:integrable_chab}
The distribution $\mu_\ab$ of the random  variable $S^\ab(t)$ can be written as
$$\mu_\ab=\kappa\,\delta_0+(1-\kappa)\tilde{\mu}_\ab$$
where $\kappa=\In_\ab\me^{-\mu t}/(\st_\ab(1-\ea\A))$ and $\tilde{\mu}_\ab$ is the distribution of a continuous random variable having probability density function
$$
\fab_S(t,x)=\frac{1}{2\pi(1-\kappa)}\int\me^{-iux}\big(\chab_t(u)-\kappa\big)\md u,
$$
with $\chab_t$  the characteristic function of $\sab(t)$ defined in (\ref{eq:charac func S(t)_A}).
\end{lemma}
\begin{proof}
 It follows from Lemma \ref{lemma:speed_thab} and Equation (\ref{eq:charac func S(t)_A}) that for all $|u|$ large
 \begin{equation}\label{eq:lem}
   \begin{split}
     \Big|\chab_t(u)-\kappa\Big|   \le&~  \Big|\me^{(\thab-\lambda-\mu)t/2}\frac{\In_\ab(\thf-\lambda-\mu)+2\lambda}{(\thab+\thf)\st_\ab(1-\ea\A)}-\frac{\In_\ab
\me^{-\mu  t}}{\st_\ab(1-\ea\A)}\Big|\\
&\quad+\Big|\me^{-(\thf+\lambda+\mu)t/2}\frac{\In_\ab(\thab+\lambda+\mu)-2\lambda}{(\thab+\thf)\st_\ab(1-\ea\A)}\Big|\\
  \le &~\Big|\me^{(\thab-\lambda-\mu)t/2}\frac{2\lambda-\In_\ab(\thab+\lambda+\mu)}{(\thab+\thf)\st_\ab(1-\ea\A)}\Big|\\
  &\quad+\Big|\frac{\In_\ab \me^{-\mu t}}{\st_\ab(1-\ea\A)}\big(\me^{(\thab-\lambda+\mu)t/2}-1\big)\Big|+C_2 \me^{-Dt u^2/2}\\
  \le &~ C_1\frac{1}{u^2}+C_2 \me^{-Dt u^2/2},
   \end{split}
 \end{equation}
 where $C_1, C_2$ are constants independent of $u$.
 This implies that the integrand in the lemma is integrable.

Without loss of generality we may suppose $\In_\ab> 0$.  Using (\ref{eq:lem}) we obtain that as $T\rightarrow\infty$
$$
\Big|\frac{1}{2T}\int_{-T}^{T}\chab_t(u)\md u-\kappa\Big|=\frac{1}{2T}\Big|\int_{-T}^{T}\Big(\chab_t(u)-\kappa\Big)   \md u\Big|\rightarrow\,0.
$$
 This implies that the point $0$ is an atom  of $\mu_\ab$, since (cf.~\cite{Billing}, page 306)
\begin{equation*}
\begin{split}
\mu_\ab(\{0\})=\lim_{T\rightarrow\infty}\frac{1}{2T}\int_{-T}^{T}\chab_t(u)\md u=\kappa.
\end{split}
\end{equation*}
Moreover, $0$ is the unique atom  of $\mu_\ab$ since (cf.~\cite{Billing}, page 306)
$$
\sum_{q}(\mu_\ab(\{q\}))^2=\lim_{T\rightarrow\infty}\frac{1}{2T}\int_{-T}^T|\chab_t(u)|^2\md u=(\mu_\ab(\{0\}))^2,
$$
where the sum is taken over the set of points of positive $\mu_\ab$ measure, and the second equality can be seen by using (\ref{eq:lem}) and the fact that $\chab_t(u)$ is uniformly bounded. This establishes  the lemma.
   \end{proof}

It follows from Lemma \ref{lemma:integrable_chab} that $\sab(t)$ is a continuous random variable if and only if $\In=(1,0)$, i.e, for instantaneous injection of the solute. It is interesting that in this case we have the following.

\begin{theorem}\label{th:1}
  The partial probability
  density functions $\hat{f}_S^\tau$  of $S^\tau(t)$ for $\tau\in\{\f,\ab\}$  satisfy the
  partial differential equations {\rm(\ref{Dequation})}:
 \begin{equation*}
 \begin{split} \frac{{\partial \jf_S(t,x)
}}{{\partial t}} +\frac{{\partial \jab_S(t,x) }}{{\partial t}}&= D\frac{{\partial^2 \jf_S(t,x)}}
{{\partial x^2 }} - v\frac{{\partial \jf_S(t,x) }}{{\partial x}},\\
\frac{{\partial \jab_S(t,x) }}{{\partial t}} &= - \mu \jab_S(t,x) +
\lambda \jf_S(t,x)
\end{split}
\end{equation*}
 for $t>0$, with initial and boundary conditions
\begin{equation*}
\begin{split}
&\jf_S(0,x):=\delta(x),\quad\jab_S(0,x):=0;\\
&\lim_{x\rightarrow\infty}\hat{f}_S^\tau(t,x)=\lim_{x\rightarrow\infty}\frac{\partial \hat{f}_S^\tau(t,x)}{\partial x}=0\quad\textrm{for}~t\ge 0, ~ \tau\in\{\f,\ab\}.
\end{split}
\end{equation*}
\end{theorem}
\begin{proof}
The initial conditions imply  $\In=(1,0)$.  It follows from Lemma \ref{lemma:integrable_chf},
\ref{lemma:integrable_chab} and Equation (\ref{eq:partial densities}) that
$$
\hat{f}^\tau_S(t,x)=\frac{1}{2\pi}\int \me^{-iux}\hat{\varphi}_t^\tau(u)\,\md u\quad\textrm{for}\quad\tau\in\{\f,\ab\},
$$
where
\begin{equation}\label{eq:partial characteristic functions}
\hat{\varphi}_t^\tau(u)=\st_\tau(1-\varepsilon_\tau\A(t))\varphi_t^\tau(u),
\end{equation}
with $\varphi^\tau_t$  the characteristic functions of $S^\tau(t)$ given in (\ref{eq:charac func S(t)_F}) and
(\ref{eq:charac func S(t)_A}) respectively. It is easy to see that $\jf_S$ and $\jab_S$ satisfy the initial and boundary conditions.

 Using Lemma
\ref{lemma:speed_thab} it is not hard to check that the four functions in $u$
$$
\Big|\frac{\partial \me^{-iux}\jchf_t(u)}{\partial t}\Big|,
\quad\Big|\frac{\partial \me^{-iux}\jchab_t(u)}{\partial t}\Big|,\quad
\Big|\frac{\partial \me^{-iux}\jchf_t(u)}{\partial x}\Big|,\quad
\Big|\frac{\partial^2 \me^{-iux}\jchf_t(u)}{\partial x^2}\Big|,
$$
are all bounded by a function of the form $C_1/u^2+C_2 \me^{-\frac{t}{2}D u^2}$ for
$|u|$ large, where $C_1, C_2$ are constants independent of $u$.
Thus we can exchange the integral and differential operators in the
partial differential equations (cf.~\cite{Durrett2010}, page 417). Hence we only need to show that
\begin{equation*}
 \begin{split} \frac{{\partial \jchf_t(u)
}}{{\partial t}}&= -Du^2\jchf_t(u)+ i v u \jchf_t(u)-
\lambda \jchf_t(u)+ \mu \jchab_t(u),\\
\frac{{\partial \jchab_t(u) }}{{\partial t}} &= - \mu \jchab_t(u) +
\lambda \jchf_t(u).
\end{split}
\end{equation*}
It follows from (\ref{eq:charac func S(t)_F}), (\ref{eq:charac func S(t)_A}), (\ref{eq:partial characteristic functions}) and $\In_\ab=0$ that
\begin{equation*}
  \begin{split}
   &-Du^2\jchf_t(u)+ i v u\jchf_t(u) -\lambda \jchf_t(u)+\mu \jchab_t(u)\\
=&~(\frac{\thab-\thf}{2}-\lambda)\Big(\me^{(\thab-\lambda-\mu)t/2}\frac{\thab-\lambda+\mu}{\thab+\thf}
+\me^{-(\thf+\lambda+\mu)t/2}\frac{\thf+\lambda-\mu}{\thab+\thf}\Big)\\
&~\hspace{1cm}+\frac{2\lambda\mu}{\thab+\thf}\big(\me^{(\thab-\lambda-\mu)t/2}-\me^{-(\thf+\lambda+\mu)t/2}\big)\\
=&~\frac{\thab-\lambda-\mu}{2}\me^{(\thab-\lambda-\mu)t/2}\frac{\thab-\lambda+\mu}{\thab+\thf}
-\frac{\thf+\lambda+\mu}{2}\me^{-(\thf+\lambda+\mu)t/2}\frac{\thf+\lambda-\mu}{\thab+\thf}\\
=&~\frac{{\partial \jchf_t(u) }}{{\partial t}},
  \end{split}
\end{equation*}
where the second equality holds since $
(\thab-\lambda+\mu)(\thf+\lambda-\mu)=4\lambda\mu. $

One finishes the proof of the theorem by checking
\begin{equation*}
\frac{{\partial \jchab_t(u) }}{{\partial t}}= - \mu \jchab_t(u) +
\lambda \jchf_t(u).
\end{equation*}

\end{proof}

We would like to point out that Lindstrom and Narasimhan \cite{Lindstrom} gave  an analytical solution of the partial differential equations with different initial and boundary conditions by using Laplace and inverse Laplace transforms. Their method can also be used with our initial and boundary conditions to give the same solutions as we have obtained via our stochastic model as   Theorem \ref{th:1}.

\section{Moments of $S(t)$ and $S^\tau(t)$}\label{sec: mom}

The mean and variance of $S(t)$ can be obtained by differentiating its characteristic function $\varphi_t$ given in
 (\ref{eq:charac func S(t)}), but a more leisurely way is to take the limits of $\expec{S_n(t)}$ and $\var{S_n(t)}$ respectively.

\begin{lemma}\label{lemma:moments convergence}
The first and second moments of $S(t)$  can be obtained by taking the limits of the corresponding moments of $S_n(t)$ respectively, i.e.,
\begin{equation*}
\begin{split}
  \expec{S_n(t)}\rightarrow\expec{S(t)},\quad \var{S_n(t)}\rightarrow\var{S(t)}
\end{split}
\end{equation*}
\end{lemma}
\begin{proof}
Recall that the mean  of $K_n$ is
given in \cite{MBD}. It is not difficult to check that the first moment of $K_n$ is uniformly bounded, i.e., there exists $M>0$,
 such that $|\expec{K_n}|\le M$.
Since the $X_k$'s are  independent random variables also independent
of $K_n$, using $\Delta t=t/n$ and (\ref{eq:displace}) we obtain
that
{\small\begin{equation*}
  \begin{split}
   \expec{S_n^3}&=\expec{\expec{S_n^3~|~K_n}} =\sum_{j=0}^n f_n(j)\expec{\Big(\sum_{k=1}^j(X_k+v\Delta t)\Big)^3}\\
   &=\expec{(X_1+v\Delta t)^3}\sum_{j=0}^n f_n(j)j+3\expec{(X_1+v\Delta t)^2}\expec{X_1+v\Delta t}\sum_{j=1}^n f_n(j)j(j-1)\\
   &\quad+\big(\expec{X_1+v\Delta t}\big)^3\sum_{j=0}^nf_n(j)j(j-1)(j-2)\\
   &\le(t+6 D vt^2+v^3 t^3)\expec{K_n}+3(2 D t+v^2 t^2)vt \expec{K_n}+v^3 t^3,
  \end{split}
\end{equation*}}
which implies that $S_n(t)$ and $S_n^2(t)$ are uniformly integrable.
This together with the fact that $S_n(t)$ converges to $S(t)$ in
distribution (shown in Section \ref{sec: towards condtious
 time})
  imply that $\expec{S_n(t)}\rightarrow\expec{S(t)}, \var{S_n(t)}\rightarrow\var{S(t)}$ (see e.g.~\cite{Billing}, Theorem 25.12).
\end{proof}

Since $X_k$ is independent of $K_n$,
from \MBDE{4} and \MBDP{2.1} together with Equation (\ref{eq:displace}) we can determine the first and second
moments of $S_n(t)$:
\begin{equation*}\label{expec_S(t)}
  \begin{split}
    \expec{S_n(t)}&=\expec{K_n}\expec{X_1+v \Delta t}=\expec{K_n}v \Delta t\\
    &=\st_\f\left(n-\frac{\ef(1-\ei^n)}{1-\ei}\right)v \Delta t=\st_\f\, v\, t-\frac{\ef\,\st_\f(1-\ei^n)}{1-\ei}v \Delta t,
  \end{split}
\end{equation*}
and
\begin{equation*}\label{var_S(t)}
  \begin{split}
   \var{S_n(t)}=&~\expec{K_n}\var{X_1+v \Delta t}+\var{K_n}(\expec{X_1+v \Delta t})^2\\
   =&~\st_\f\Big(n-\frac{\ef(1-\ei^n)}{1-\ei}\Big) 2 D\Delta t+\var{K_n}(v\Delta t)^2\\
   =&~2 D\st_\f\, t-\frac{2 D\ef\st_\f(1-\ei^n)}{1-\ei}\Delta t
   +\frac{\st_\ab(1+\ei)+2\ef(\st_\ab-\st_\f)\ei^n}{1-\ei}\st_\f\,v^2 t\Delta t\\
   &~\hspace{3cm}+\frac{\ei\big(\ef(\st_\f-\st_\ab)-2\st_\ab\big)-\ef(\st_\ab-\In_\f)}{(1-\ei)^2}\st_\f(v \Delta t)^2\\
   &~\quad+\ei^n\Big(\frac{\ef(\st_\f-\st_\ab)}{1-\ei}+2\frac{\ei\st_\ab+\ef(\st_\ab-\In_\f)}{(1-\ei)^2}
   -\ei^n\frac{\st_\f \ef^2}{(1-\ei)^2}\Big)\st_\f(v \Delta t)^2.
  \end{split}
\end{equation*}
Substituting $\Delta t=t/n$ and (\ref{Delta t}) in the mean and variance of $S_n(t)$ and letting $n\rightarrow \infty$, by Lemma \ref{lemma:moments convergence} we
obtain  the mean and variance of $S(t)$.
\begin{proposition}\label{prop:mean var S(t)}
The mean and variance of $S(t)$ are given by
\begin{equation*}\label{expec_S(t)_lambda}
\expec{S(t)}=\st_\f v\, t-\frac{\ef\st_\f}{\lambda+\mu}v (1-\A),
\end{equation*}
and
\begin{equation*}\label{var_S(t)_lambda}
\begin{split}
  \var{S(t)}=& \;2 D \st_\f t-\frac{2 D \ef \st_\f}{\lambda+\mu}(1-\A)+\frac{2\big(\st_\ab+\ef(\st_\ab-\st_\f)\A\big)}{\lambda+\mu}\st_\f v^2 t\\
&\quad+\frac{\ef(\st_\f-\st_\ab)-2\st_\ab-\ef(\st_\ab-\In_\f)}{(\lambda+\mu)^2}\st_\f v^2\\
&\quad\qquad+\A\Big(2\frac{\st_\ab+\ef(\st_\ab-\In_\f)}{(\lambda+\mu)^2}-\A\frac{\st_\f \ef^2}{(\lambda+\mu)^2}\Big)\st_\f v^2.
\end{split}
\end{equation*}

\end{proposition}

 Now we are going to consider the means and variances of $S^\tau(t), \tau\in\{\f,\ab\}.$ Again, one could obtain them from their characteristic functions, but  we will use the following lemma, which can be proved in a similar way as Lemma \ref{lemma:moments convergence}.
\begin{lemma}\label{lemma:conditional moments convergence}
The first and second moments of $S^\tau(t), \tau\in\{\f,\ab\}$ can be obtained by taking the limits of the corresponding moments of  $S^\tau_n(t)$, i.e.,
\begin{equation*}
\begin{split}
 \expec{S_n^\tau(t)}\rightarrow\expec{S^\tau(t)},\quad \var{S_n^\tau(t)}\rightarrow\var{S^\tau(t)}.
\end{split}
\end{equation*}
\end{lemma}

Because of independence, using  \MBDE{5} and
\MBDP{3.1} together with Equation (\ref{eq:displace}), we obtain that
\begin{equation*}\label{expecfr_S(t)}
\expec{S^\f_n(t)}=\expec{{K}_n^\f}v\Delta t
=~\frac{\st_\f-\ef\st_\ab\,\ei^{n-1}}{1-\ef\ei^{n-1}} v t+\frac{(\st_\ab-\ef\st_\f
)(1-\ei^n)}{(1-\ei)(1-\ef\ei^{n-1})} v \Delta t,
\end{equation*}
and
\begin{equation*}\label{varfr_S(t)}
\begin{split}
&\var{S^\f_n(t)} =   \expec{{K}_n^\f} \var{X_k+v\Delta t}+\var{{K}_n^\f}(\expec{X_k+v\Delta t})^2\\
           \quad\quad=&   \frac{\st_\f-\ef\st_\ab\,\ei^{n-1}}{1-\ef\ei^{n-1}} 2D t
          +\frac{(\st_\ab-\ef\st_\f \, )(1-\ei^n)}{(1-\ei)(1-\ef\ei^{n-1})} 2D \Delta t
          +\frac{\st_\f^2-\ef\st_\ab^2\,\ei^{n-1}}{1-\ef\ei^{n-1}}v^2 t^2\\
          &~\hspace{4cm}-\left(\frac{\st_\f-\ef\st_\ab\,\ei^{n-1}}{1-\ef\ei^{n-1}}\,t
  +\frac{(\st_\ab-\ef\,\st_\f)(1-\ei^n)}{(1-\ei)(1-\ef\ei^{n-1})}\,\Delta t\right)^2\! v^2\\
  &\hspace{1cm}-\left(\frac{\st_\ab\,\st_\f(1+3\ef\ei^{n-1})}{1-\ef\ei^{n-1}}
  +2\,\frac{\ef\,\st_\f^2+\st_\ab^2\ei^n-2\st_\ab\,\st_\f(1+\ef\ei^{n-1})}{(1-\ei)(1-\ef\ei^{n-1})}\right)v^2 t\,\Delta t \\
  +&~(1-\ei^n)\left(\frac{\st_\ab\st_\f(4+\ef)-(\st_\ab+\ef\st_\f^2)}{(1-\ei)(1-\ef\ei^{n-1})}
  +2\,\frac{\ef\,\st_\f^2+\st_\ab^2-2\st_\ab\,\st_\f(1+\ef)}{(1-\ei)^2(1-\ef\ei^{n-1})}\right)(v\Delta t)^2.
\end{split}
\end{equation*}
Substituting $\Delta t=t/n$ and (\ref{Delta t}) in the mean and variance of $S^\f_n(t)$ and letting $n\rightarrow \infty$, by Lemma \ref{lemma:conditional moments convergence} we
obtain the mean and variance of $S^\f(t)$.
\begin{proposition}\label{prop:mean var SF(t)}
The mean and variance of $S^\f(t)$ are given by
\begin{equation*}\label{expecfr_S(t)_lambda}
\expec{S^\f(t)}=\frac{\st_\f-\ef\,\st_\ab \A}{1-\ef \A} v t+\frac{(\st_\ab-\ef\,\st_\f)(1-\A)}{(\lambda+\mu)(1-\ef \A)}v,
\end{equation*}
and
\begin{equation*}\label{varfr_S(t)_lambda}
\begin{split}
&\var{S^\f(t)}  =\frac{\st_\f-\ef\st_\ab \A}{1-\ef \A}2 D t+\frac{\st_\ab-\ef\st_\f}{(\lambda+\mu)(1-\ef \A)}2D(1-\A)+\frac{\st_\f^2-\ef\st_\ab^2 \A}{1-\ef
\A}v^2 t^2 \\
&\hspace{5cm}-\Big(\frac{\st_\f-\ef\st_\ab \A}{1-\ef
\A}t+\frac{\st_\ab- \ef\st_\f}{(\lambda+\mu)(1-\ef
\A)}(1-\A)\Big)^2 v^2\\
&\quad-2\frac{\ef\st_\f^2+\st_\ab^2 \A-2\st_\ab\st_\f(1+\ef
\A)}{(\lambda+\mu)(1-\ef \A)}v^2
t+2(1-\A)\frac{\ef\st_\f^2+\st_\ab^2-2\st_\ab\st_\f(1+\ef)}{(\lambda+\mu)^2(1-\ef
\A)}v^2.
\end{split}
\end{equation*}

\end{proposition}

In a quite similar way we obtain the following result.
\begin{proposition}\label{prop:mean var SA(t)}
The mean and variance of $\sab(t)$ are given by
\begin{equation*}\label{expecab_S(t)_lambda}
\expec{S^\ab(t)}=\frac{\st_\f-\ea\,\st_\ab \A}{1-\ea \A} v t+\frac{(\ea\st_\ab-\st_\f)(1-\A)}{(\lambda+\mu)(1-\ea \A)}v,
\end{equation*}
and
\begin{equation*}\label{varab_S(t)_lambda}
\begin{split}
&\var{S^\ab(t)}  =\frac{\st_\f-\ea\st_\ab \A}{1-\ea A}2 D t+\frac{\ea\st_\ab-\st_\f}{(\lambda+\mu)(1-\ea \A)}2D(1-\A)+\frac{\st_\f^2-\ea\st_\ab^2 \A}{1-\ea \A}v^2 t^2 \\
&\hspace{5.2cm}-\Big(\frac{\st_\f-\ea\st_\ab \A}{1-\ea \A}t+\frac{\ea\st_\ab-\st_\f }{(\lambda+\mu)(1-\ea \A)}(1-\A)\Big)^2 v^2\\
&\;-2\frac{\st_\f^2+\ea\st_\ab^2 \A-\st_\ab\st_\f(1+\ea)(1+\A)}{(\lambda+\mu)(1-\ea \A)}v^2
t+2(1-\A)\frac{\st_\f^2+\ea\st_\ab^2-2\st_\ab\st_\f(1+\ea)}{(\lambda+\mu)^2(1-\ea
\A)}v^2.
\end{split}
\end{equation*}

\end{proposition}

Now, we will use our model to illustrate some mistake made by
Michalak and Kitanidis \cite{Mich}.
Note that the moments of $S^\f(t)$ and $S^\ab(t)$ calculated in Proposition \ref{prop:mean var SF(t)} and \ref{prop:mean var SA(t)} are exactly the same as that calculated by Michalak and Kitanidis (see Section \ref{sec: PDE}) if the initial conditions are specified as $\In=(1,0)$, since for $m\ge 1$
$$
\mathrm{E}_{(1,0)}[(S^\tau(t))^m]=\int x^m{f_S^\tau(t,x)}\md
x=\int\limits{x^{m}\frac{C_\tau(t,x)}{{M}_{\tau}^{ ( 0)
}(t)}
\,\mathrm{d}x}={M}_{\tau}^{(m)}(t),\quad\tau\in\{\f,\ab\}
$$
where the second equality holds since by Theorem \ref{th:1} both the  partial probability density functions $\hat{f}^\tau_S=\st_\tau(1-\varepsilon_\tau\A)f_S^\tau$ and the concentration functions $C_\tau$ satisfy the partial differential equations (\ref{Dequation}) with the same initial and boundary conditions.

 Recall from Section \ref{sec:
PDE} that we translate the parameters into our paper as follows:
$$
 \mu=k,\quad\lambda=\beta k.
 $$
 If we let the solute  be
`\free' at time $0$ and $t$, i.e., the initial distribution
$\In=(1,0)$, then
$$\ef=-\frac{\lambda}{\mu}=-\beta, \quad\st_\f=\frac{\mu}{\lambda+\mu}=\frac{1}{\beta+1},\quad \st_\ab=\frac{\lambda}{\lambda+\mu}=\frac{\beta}{\beta+1}.$$
 Substituting these parameters into Proposition \ref{prop:mean var SF(t)} yields
\begin{equation*}
\begin{split}
{\rm Var}_{(1,0)}(S^\f(t))=&\;\frac{1+\beta^2 \A}{(\beta+1)(1+\beta \A)}2 Dt+\frac{2\beta}{k(\beta+1)^2(1+\beta \A)}2D(1-\A) \\
 +&\frac{1+\beta^3 \A}{(\beta+1)^2(1+\beta \A)}v^2 t^2-\Big(\frac{1+\beta^2\A}{(\beta+1)(1+\beta \A)}t+\frac{2\beta(1-\A)}{k(\beta+1)^2(1+\beta
 \A)}\Big)^2v^2\\
 &\hspace{2.0cm}-\frac{6\beta(\beta \A-1)}{k(\beta+1)^3(1+\beta \A)}v^2 t+\frac{6\beta(\beta-1)}{k^2(\beta+1)^4(1+\beta
 \A)}v^2(1-\A),
\end{split}
\end{equation*}
where $\A=\exp{(-(\lambda+\mu)t)}=\exp{(-(\beta+1)k t)}$.  This gives indeed Equation~(\ref{eq:Michalak}) which is taken from \cite{Mich}.

However, Michalak and Kitanidis state in their paper that
$\var{S^\tau(t)}$ can be obtained by a linear combination of
$\varfr{S^\tau(t)}$ and $\varab{S^\tau(t)}$ (i.e., $\var{S^\tau(t)}$ with
initial distributions $\In=(1,0)$ and $\In=(0,1)$). This is not
true, and we provide the correct formulas in Proposition
\ref{prop:mean var SF(t)} and \ref{prop:mean var SA(t)}. We have also
given the formula for the total solute in Proposition \ref{prop:mean
var S(t)}.

\section{Double-peak behavior in reactive transport models}\label{sec:Double Peaking}
Double peaks in the `free' concentration distribution $C_\f$ are discussed by
Michalak and Kitanidis \cite{Mich} using simulations. Theorem \ref{th:1} tells us that $C_\f(t,\cdot)$ can be seen as the partial probability density function $\jf_S(t,\cdot)$ of $S^\f(t)$ if the initial distribution is $\In=(1,0)$. We will show in this
section how double peaks can also be explained by means of
our stochastic reactive transport model. Let
$ \ff_{S_n}(t,\cdot)$ be the probability density function of $\sf_n(t)$ defined in Section \ref{sec: SRTM}. We are going to approximate $\ff_S(t,\cdot)$ by $\ff_{S_n}(t,\cdot)$, since $\sf_n(t)$ converges to $\sf(t)$ in distribution.

Michalak and Kitanidis consider Gaussian diffusion, i.e., the
$X_k$'s are normally distributed random variables with mean $0$ and
variance $2 D \Delta t$, which satisfy Equation (\ref{eq:displace}).
So the
characteristic function
 of $\sf_n(t)$ can be written  as
\begin{equation*}
  \chf_{t,n}(u)=\expec{\me^{i u {\sf_n(t)}}}=\gf_n\Big(\expec{\me^{i u(X_1+ v\Delta t)}}\Big)=\sum_{j=0}^n \ff_{n}(j)\exp{\big(i v \Delta tj u-   D \Delta tj u^2\big)},
\end{equation*}
where $\ff_n$ is the probability mass function of $K_n^\f$
defined in Equation (\ref{eq:fntau}). Obviously $\int_{-\infty}^{\infty}| \chf_{t,n}(u)|\mathrm{d}u<\infty$.
Thus by the inverse Fourier transformation, using that $\ff_n(0)=0$, we obtain
\begin{equation}\label{f_{S_n}}
\begin{split}
\ff_{S_{n}}(t,x)=&~\frac{1}{2\pi}\int_{-\infty}^{\infty}\me^{-i u
x}\chf_{t,n}(u)\mathrm{d}u\\
=&~\sum_{j=0}^n \ff_n(j) \frac{1}{2\pi}\int_{-\infty}^{\infty}\exp{\Big(i u (j v \Delta t-  x )- u^2 j D  \Delta t \Big)}\mathrm{d}u\\
=&~\sum_{j=1}^n  \ff_n(j)\frac{1}{\sqrt{  4 \pi j D \Delta t
}}\exp{\Big(-\frac{(x-j v  \Delta t)^2}{4  j D \Delta t}\Big)}.
\end{split}
\end{equation}
 So $\sf_n(t)$ is a mixture of Gaussian distributions with mean $j v\Delta t$ and variance $2 j D \Delta t$. Recall from \cite{MBD} that the probability mass function $\ff_n$ of $K_n^\f$ can be unimodal or bimodal. This property of $K_n^\f$ gives rise to the same phenomenon for $\sf_n(t)$, i.e., one peak or two peaks appear in the probability density function $\ff_{S_n}(t,x)$ of $\sf_n=\sf_n(t)$ for large $n$.
\begin{figure}[h]
  \centering{ \includegraphics[width=3.5cm]{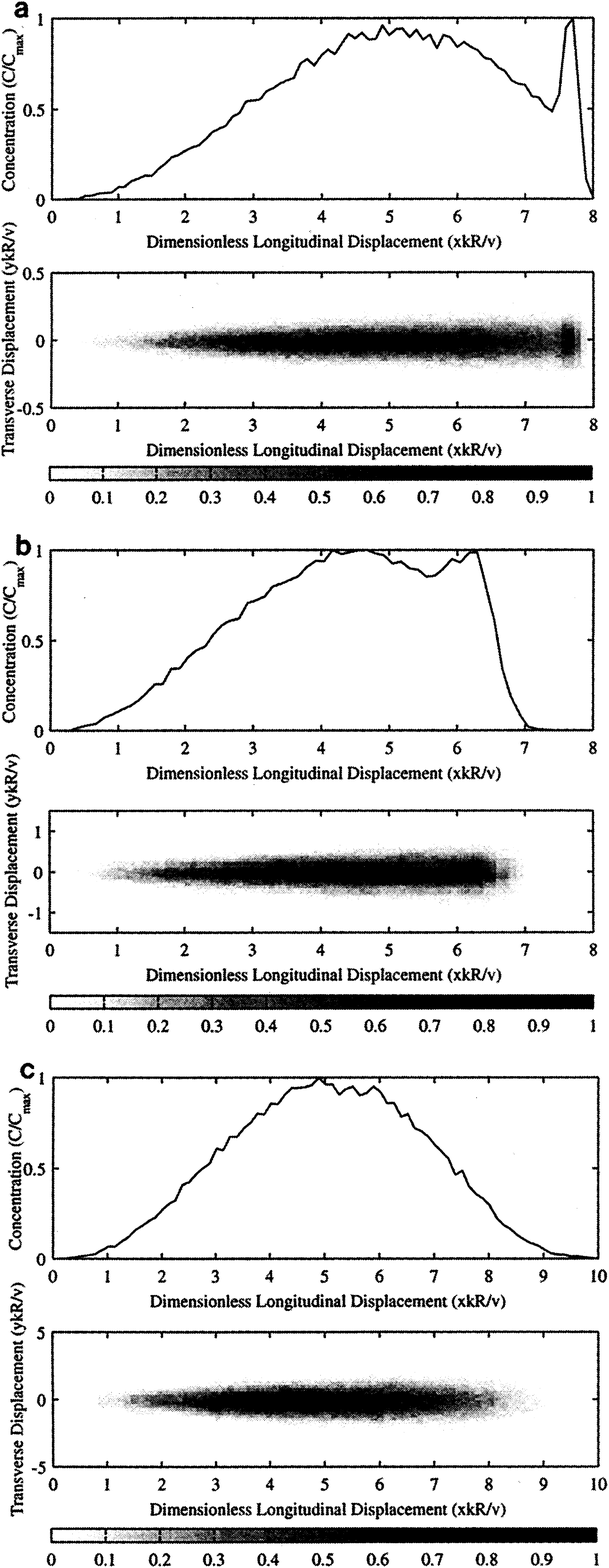}\quad\quad\quad\includegraphics[width=5.0cm]{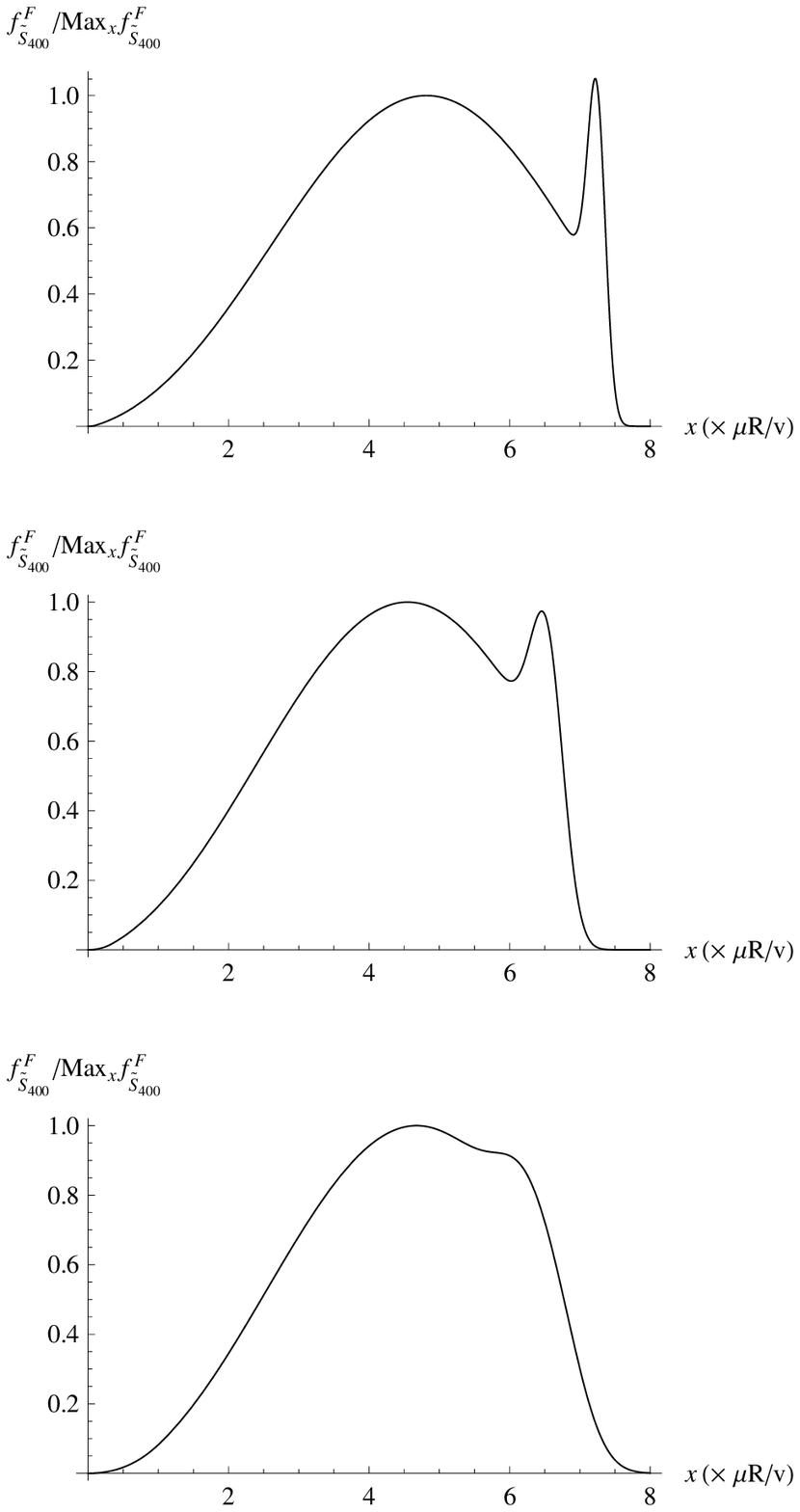}\\
  \caption{The three graphs in the left column are the normalized concentration functions $C_\f(t,\cdot)/\max_{x}C_\f(t, x)$ copied from Michalak
  and Kitanidis \cite{Mich}. The three graphs in the right column are the normalized probability density functions $f^\f_{\tilde{S}_{400}}(t,\cdot)/max_x f^\f_{\tilde{S}_{400}}(t,x)$ given by the Fourier transformation in our paper. All graphs have $Pe=100, v=L=1,~ R=2$.   In the first row $ Da_I= 0.1, t^*=3.6$, in the second row $ Da_I= 0.33, t^*=3.2$, and in the last row $ Da_I= 1.0, t^*=3.0$.}\label{pic10}}
\end{figure}

Michalak and Kitanidis  focus on the case
that the solute starts in the free phase and the length of the initial solute is $L$, i.e., the initial distributions of the PDE's (\ref{Dequation}) are given by
$$
C_\f(0,x)=\frac{1}{L}\textbf{1}_{[0,L]}(x), \quad C_\ab(0,x)=0.
$$
So to make the comparison, we look at the probability density function $f_{\tilde{S}_n}^\f(t,x)$ of
$$
\tilde{S}^\f_n(t)=S^\f_n(t)+U_L,
$$
where $U_L$ is a uniformly distributed random variable over the interval $[0,L]$ (independent of $S_n^\f(t)$). Michalak and Kitanidis point out that the double peaking behavior of the free concentration distribution is a function of the so called
 Damk\"{o}hler number
of the first kind $ Da_I=\mu L R/v, $ where $R$ is the dimensionless retardation
coefficient. They state that the timing of its appearance is
controlled by the mass transfer rate and the retardation factor,
i.e., the dimensionless time $t^*=\mu(R-1)t$. The so called P\'{e}clet
number $Pe=v L/D $ is  kept constant at a value of $100$.
Recalling from Section \ref{sec: PDE} that $\lambda=\beta \mu=(R-1)\mu$ and $a=\lambda\Delta t, b=\mu\Delta t$, we translate these
parameters into our paper as follows:
$$
a=\frac{k(R-1)t}{n}=\frac{t^*}{n},\quad b=\frac{k\, t}{n}=\frac{t^*}{n(R-1)},\quad D=\frac{v L}{Pe},\quad
\Delta t=\frac{t}{n}=\frac{t^*L R}{ nv(R-1)Da_I }.
$$

The graphs in the left column of Figure \ref{pic10} are a copy of
the graphs of the normalized aqueous concentration functions $C_\f(t,\cdot)/\max_{x}C_\f(t, x)$ (consisting of the free particles) in
Michalak and Kitanidis \cite{Mich} using simulations corresponding
to different choices of the Damk\"{o}hler number $Da_I$ and
dimensionless time $t^*$. The  three graphs in the right column of
Figure \ref{pic10} are the normalized  density functions $f^\f_{\tilde{S}_{400}}(t,\cdot)/\max_x f^\f_{\tilde{S}_{400}}(t, x)$
 calculated using Equation (\ref{f_{S_n}}) corresponding
to the same choice of $Da_I$ and $t^*$. The number $n$ is chosen large enough such that $\max{(a, b)}=\max{(\lambda \Delta t,\mu\Delta t)}\le 0.01$. From Figure \ref{pic10} it
is obvious that our model gives a much better view at the double
peaking phenomenon.
\begin{figure}[h]
  \centering {\includegraphics[width=7.0cm]{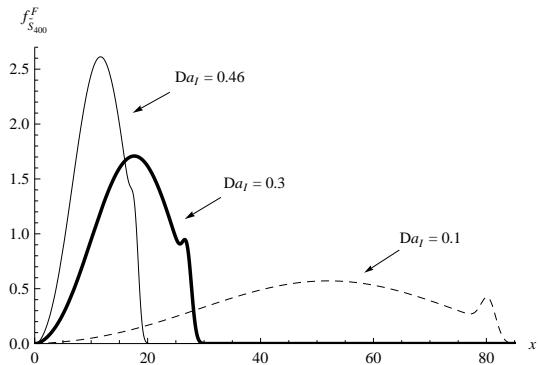}\\
  \vspace{0cm}
  \caption{The three graphs  are the probability density functions $\ff_{\tilde{S}_{400}}(t,\cdot)$ of $\tilde{S}_{400}^\f(t)$. All graphs have $Pe=100, v=L=1, R=2, t^*=4.0$, and different Damk\"{o}hler numbers.\label{pic12}}}
\end{figure}

Moreover, for each $t^*$, by a numerical calculation we can obtain
upper bounds for $Da_I$ such that double peaks appear. For  example, Figure \ref{pic12} gives an intuition on how double peaks behave when $Da_I$ increases. We numerically calculated the upper bounds for $Da_I$ in Table $1$  corresponding to different
dimensionless times $t^*$ with $R=2$. For example, when $t^*=2.0$
two peaks occur for all $Da_I>0$ until $Da_I= Da_I^{\max}=0.43$.
Table $1$ suggests that double peaking is pronounced for $2\le
t^*\le 5$, and almost dies out when $t^*<1.5$  or $t^*>10$.

\begin{center}
Table 1. $R=2, Pe=100, v=1, L=1, n=400.$

{\small\begin{tabular}{c|c|c|c|c|c|c|c|c|c|c|c|c|c}
  \hline
  $t^*$&1.5&2.0&2.5&3.0&3.5&4.0&4.5&5.0&6.0&7.0&8.0&9.0&10.0\\\hline
  $Da_I^{\max}$&0.12&0.43&1.45&1.42&0.73&0.45&0.30&0.21&0.11&0.07&0.04&0.02&0.02\\\hline
\end{tabular}}
\end{center}

\section{Final remarks}

 We emphasize that the so called `random walk method'
or `particle tracking method' first proposed by Kinzelbach \cite
{Kinzelbach88} has a relation to our model, but has always been used
as a simulation tool, to perform numerical experiments (for a recent
example see \cite{Meer}). In fact it is shown in \cite{AKAD} for the
first time that if one takes an appropriate limit (in a similar way as
in \cite{dehling}), then the Fokker-Planck equations of an extended version of our simple model to a Markov chain which also involves discrete steps in space, yield the partial differential
equations (\ref{Dequation}) in Section~\ref{sec: PDE}.

Finally we mention that our computations yield the following. If one starts in  the stationary distribution, i.e., $\In=(\st_\f,\st_\ab)$, then $\ef=\ea=0$. Substituting
$$\ef=0,\quad \st_\f=\frac{\mu}{\lambda+\mu},\quad \st_\ab=\frac{\lambda}{\lambda+\mu}$$
in Proposition \ref{prop:mean var S(t)}, we obtain
\begin{equation*}
  \begin{split}
 \var{S(t)}=&\frac{2 D \mu}{\lambda+\mu}t+\frac{2\mu\lambda}{(\lambda+\mu)^3}v^2 t-\frac{2\mu\lambda}{(\lambda+\mu)^4}v^2\Big(1-\exp{\big(-(\lambda+\mu)t\big)}\Big).
  \end{split}
\end{equation*}
We then recuperate a (more general and more detailed) version of
the main result of Gut and Ahlberg (\cite{gut}, p.251).

\bibliography{Markov-bin-AKAD}
\bibliographystyle{plain}
\end{document}